\documentclass[10pt,journal,twocolumn,final]{IEEEtran}

\usepackage{amssymb,amsmath,mathpazo}
\usepackage{graphicx,framed}
\usepackage[caption=false]{subfig}
\usepackage{algorithmic}

\title{Enhanced image approximation using shifted rank-$1$ reconstruction}
\author{F. Bo\ss mann\thanks{Harbin Institute of Technology, China ({florian.bossmann@uni-passau.de})},
   \and J. Ma\thanks{Harbin Institute of Technology, China ({jma@hit.edu.cn})}}

\newcommand{\C}{\mathbb{C}}
\newcommand{\Z}{\mathbb{Z}}

\newtheorem{theorem}{Theorem}
\newtheorem{definition}{Definition}
\newtheorem{remark}{Remark}
\newtheorem{algorithm}{Algorithm}

\begin{document}

\maketitle

\begin{abstract}
 Low rank approximation has been extensively studied in the past. It is most suitable to reproduce rectangular like structures in the data. In this work we introduce a generalization using "shifted" rank-$1$ matrices to approximate $A\in\C^{M\times N}$. These matrices are of the form $S_{\lambda}(uv^*)$ where $u\in\C^M$, $v\in\C^N$ and $\lambda\in\Z^N$. The operator $S_{\lambda}$ circularly shifts the $k$-th column of $uv^*$ by $\lambda_k$. These kind of shifts naturally appear in applications, where an object $u$ is observed in $N$ measurements at different positions indicated by the shift $\lambda$. The vector $v$ gives the observation intensity. Exemplary, a seismic wave can be recorded at $N$ sensors with different time of arrival $\lambda$; Or a car moves through a video changing its position in every frame. We present theoretical results as well as an efficient algorithm to calculate a shifted rank-$1$ approximation in $O(NM\log M)$. The benefit of the proposed method is demonstrated in numerical experiments. A comparison to other sparse approximation methods is given. Finally, we illustrate the utility of the extracted parameters for direct information extraction in several applications including video processing or non-destructive testing.
\end{abstract}

\begin{IEEEkeywords}
  low rank approximation, singular value decomposition, shift-invariant dictionary learning, column shifts, adaptive approximation
\end{IEEEkeywords}

\section{Introduction}

Low rank approximation of a data matrix $A\in\C^{M\times N}$ is of great interest in many applications. It is used e.g., for video denoising \cite{video}, object detection \cite{object}, subspace segmentation \cite{subspace} or seismic data processing \cite{seismic1, seismic2}. The problem can be formulated as follows: Given the matrix $A$ and a small number $L$ one seeks for vectors $u^1,\ldots,u^L\in\C^M$, $v^1,\ldots,v^L\in\C^N$ and coefficients $\sigma_1,\ldots,\sigma_L$ such that
\begin{align}\label{lowRankApprox}
\|A-\sum_{k=1}^L\sigma_ku^k(v^k)^*\|_F
\end{align}
is small. Here $\|\cdot\|_F$ is the Frobenius norm. The minimum of (\ref{lowRankApprox}) is given by the singular value decomposition  (SVD) \cite{SVD} of $A=U\Sigma V^*$ where $U\in\C^{M\times M}$, $V\in\C^{N\times N}$ are orthogonal matrices and $\Sigma\in\C^{M\times N}$ is a diagonal matrix. Then the best rank $L$ approximation is given by the first $L$ columns $u^1,\ldots,u^L$ of $U$, the first $L$ columns of $V$ and the corresponding singular values $\sigma_1,\ldots,\sigma_L$. In cases where calculating the SVD is computational too costly, randomized algorithms such as \cite{fastSVD1,fastSVD2} may be used.

In some applications only a few coefficients or affine samples of the matrix $A$ may be given. Recovering its low-rank structure is referred to as low-rank matrix completion or matrix sensing. Several techniques have been developed to solve this problem, e.g., iterative thresholding \cite{MC1}, alternating minimization \cite{AlterMin}, nuclear norm minimization \cite{undersampledSVD,NNLR} or procrustes flow \cite{ProFlo}. Theoretical results can also be found in  \cite{MRls,TheoNNM}. An overview of modern techniques is given in \cite{MRs}.

Other works consider a modified version of the low-rank approximation problem. Exemplary, a weighted Frobenius norm \cite{weightedLR}, simultaneous approximation of several matrices \cite{multiLR} or structure preserving techniques \cite{structuredLR1,structuredLR2} have been analyzed. More complex models such as block low-rank \cite{BLR}, non-negative matrix factorization \cite{NMF} or arbitrary data types \cite{GLR} have also been considered. Low-rank tensor approximation is discussed e.g., in \cite{tensorNP,TLR1,TLRs}.

Recently, ideas to combine low-rank and sparse approximation came up. The "low-rank plus sparse" model assumes that the given matrix can be written as sum of a low-rank and a sparse matrix, i.e., a matrix with only a few non-zero entries. This model is suitable e.g., in MRI \cite{LRS1}. Several techniques have been developed to deal with this problem such as \cite{LRS2,LRS5,LRS6,LRS4}. Theoretical limitations are discussed in \cite{LRS3}.

In this work, we present a generalization of low rank approximation that can be seen as "low-rank plus shift". We seek to find an approximation of $A$ using $L$ rank-$1$ matrices $u^k(v^k)^*$ where the columns of each matrix can be shifted arbitrarily. Let $\lambda^k\in\Z^N$  and $S_{\lambda^k}:\C^{M\times N}\rightarrow\C^{M\times N}$ be the operator that circularly shifts the $j$-th column of a matrix by $\lambda^k_j$. In this work we consider the following problem
\begin{align}\label{multishift}
\min\limits_{\substack{\lambda^1,\ldots,\lambda^L, \\ u^1,\ldots,u^L \\ v^1,\ldots,v^L}}
\|A-\sum\limits_{k=1}^LS_{\lambda^k}(u^k(v^k)^*)\|_F,
\end{align}
where we call $S_{\lambda^k}(u^k(v^k)^*)$ a shifted rank-$1$ matrix. The motivation behind this new approach is given by applications such as seismic data processing \cite{seismic1, seismic2}, where the given data $A$ consists of several time signals recorded at different locations. Assume each time signal records the same events (such as earthquakes or explosions) where each event is shifted according to the distance between source and sensor. Then an event can be written as such a shifted rank-$1$ matrix where $u^k$ is the seismic wave, $v^k$ is the intensity and $\lambda^k$ gives the time of arrival at the different sensors for each event. This model also holds for other applications, exemplary shown in the numerical experiments for ultrasonic non-destructive testing and video processing. We also give a more extensive motivation in Section \ref{Sec::Algorithm} after the details of the algorithm have been clarified.

Problem (\ref{multishift}) contains low-rank approximation as the special case where $\lambda^k=0$ for all $k\leq L$. Furthermore, if we assume that $\lambda^1=\lambda^2=\ldots=\lambda^k=:\lambda$, then  (\ref{multishift}) simplifies to minimizing $\|A-S_{\lambda}(UV^*)\|_F$ where $U\in\Z^{M\times L}$ and $V\in\Z^{L\times N}$. This formulation looks quite similar to the above mentioned matrix sensing problem where $A$ is the sampling of a rank $L$ matrix using the linear operator $S_{\lambda}$. However, note that in our model the shift vector $\lambda$ and thus the linear operator $S_{\lambda}$ is unknown while in matrix sensing it is assumed that the affine operator is known. Moreover, we are especially interested in the case where $\lambda^k\neq\lambda^j$ for $k\neq j$, i.e., where each rank-$1$ matrix can be shifted individually.

Indeed, problem (\ref{multishift}) is more closely related to shift-invariant dictionary learning. The dictionary learning problem can be formulated as follows. Given data $A$ find a dictionary $D=[d^1,\ldots,d^R]$ such that $A=DX$ where $d^k\in\C^M$ and $X$ is a matrix with $L$-sparse columns. The number $R$ of dictionary elements (so called atoms) can be chosen arbitrary depending on the application. Shift-invariant learning seeks for a dictionary where for each atom $d^k$ also all of its circulant shifts are contained in $D$. Let $D=[D^1,\ldots,D^R]$ where $D^k$ is a matrix with $d^k$ and all its shifts as columns. Accordingly set $X^T=[(X^1)^T,\ldots,(X^R)^T]$. Then
\begin{align*}
A=DX=\sum\limits_{k=1}^R D^kX^k.
\end{align*}
Now $D^kX^k$ can be written as shifted rank-$1$ matrix $S_{\lambda^k}(d^k(x^k)^*)$ if and only if each column of $X^k$ is at most $1$-sparse. The non-zero elements of each column combined form the vector $x^k\in\C^N$, the shift $\lambda^k$ is given by the row index of the non-zero elements. Thus (\ref{multishift}) can be formulated as shift-invariant dictionary learning where $R=L$ and each atom can only be chosen once per column of $A$. This uniqueness of each atom will play an important role in our algorithm since it allows the unique identification of events in the data and thus direct segmentation or tracking of objects.

An overview of the methods and applications in dictionary learning can be found in \cite{DLoverview}. As algorithms that solve the shift-invariant dictionary learning problem we exemplary mention the union of circulants dictionary learning algorithm (UC-DLA) \cite{UCDLA} and the matching of time invariant filters (MoTIF) \cite{MOTIF}. MoTIF aims to construct a set of uncorrelated filters. We will also use this algorithms as comparison in our numerical experiments.

The remainder of this work is structured as follows. In the next section we define the used mathematical notation and introduce the shift operators. Some important basic properties of these operators are shown. Afterwards in Section \ref{Sec::Algorithm} the idea and workflow of the proposed greedy method will be explained. Section \ref{Sec::Greedy} discusses the technical details of the greedy choice. Here we also analyze how the algorithm can be implemented efficiently in Fourier domain. Afterwards a detailed numerical evaluation is made, including the applicability in practical problems. We end our work with a conclusion and discuss how the method can be improved in future works.

\section{Preliminaries}\label{Sec::Preliminaries}

In this paper data is assumed to be given as a matrix $A\in\C^{M\times N}$ where $M,N\in\mathbb{N}$. We set $K=\min(M,N)$. The singular values of matrices will be denoted by $\sigma_k$ with $k=1,\ldots,K$, where it is stated in the text to which matrix we refer. We use vectors $u\in\C^M$, $v\in\C^N$ and $\lambda\in\Z^N$ to construct our shifted rank-$1$ matrices. In case of dealing with several vectors of the same type, we enumerate them using upper indices $k=1,\ldots,L$, i.e., we use $u^k$, $v^k$ and $\lambda^k$. Here $L$ denotes the number of shifted rank-$1$ matrices used for the approximation. We refer to element $j$ of a vector by using a lower index, e.g., $u_j$ or $u^k_j$ refer to the $j$-th index of the vector $u$ respectively $u^k$. For matrices, which are denoted by upper case letters, we use lower case to refer to their entries, e.g., $A=[a_{jk}]_{j,k=1}^{M,N}$. Square brackets indicate a matrix definition.

Throughout the work, the following matrix- or vector-operations are used: Complex conjugate ($\overline{A}$,$\overline{u}$), transpose ($A^T$,$u^T$), conjugate transpose ($A^*$,$u^*$), element wise absolute value ($|A|,|u|$), inverse of a matrix or operator ($A^{-1}$), rank of a matrix ($\mathop{rank}(A)$). We use braces whenever this notation conflicts with an upper index, e.g., $(u^k)^T$ is the transpose of $u^k$.

We consider the Frobenius norm $\|\cdot\|_F$ and the spectral norm $\|\cdot\|_2$ for matrices as well as the euclidean norm $\|\cdot\|_2$ for vectors, which are defined as
\begin{align*}
\|A\|_F^2=\sum\limits_{j,k=1}^{M,N}|a_{jk}|^2=\sum\limits_{k=1}^{K}\sigma_k^2, &&
\|A\|_2^2=\sigma_1^2,\\
\|u\|_2^2=\sum\limits_{j=1}^M|u_j|^2.
\end{align*}
The spectral norm and the euclidean norm share the same index because the matrix norm is induced by the vector norm. Let $P\in\C^{M\times N}$ be another matrix. We define the according inner product for matrices and vectors as
\begin{align*}
\langle A,P\rangle = \mathop{Trace}(A^*P), &&
\langle u^1,u^2\rangle = (u^1)^*u^2.
\end{align*}
Here $\mathop{Trace}$ is the sum of diagonal elements. Furthermore, we notate the element-wise (Hadamard-)product of two matrices by $A\odot P$.

Now we can introduce the shift operator $S_\lambda$. Therefore, let $\mathbb{I}\in\C^{(M-1)\times(M-1)}$ be the identity matrix and define
\begin{align*}
\tilde{S}:=\begin{bmatrix}
0 & 1 \\ \mathbb{I} & 0
\end{bmatrix}\in\mathbb{C}^{M\times M},&&
\tilde{S}^k:=\prod\limits_{j=1}^k S.
\end{align*}
Then $\tilde{S}u$ forward shifts the elements of $u$ by one and $\tilde{S}^k$ respectively shifts it by $k$ entries. For $k\leq0$ we use the identity $\tilde{S}^k=\tilde{S}^{M+k}$.

\begin{definition}[Shift operators]\label{defShiftOP}
Denote the columns of the matrix $A$ by $A=[a^1,\ldots,a^N]$. For $\lambda\in\Z^N$ define the shift operator $S_{\lambda}:\C^{M\times N}\rightarrow\C^{M\times N}$ as
\begin{align*}
S_{\lambda}A:=[\tilde{S}^{\lambda_1}a^1,\ldots,\tilde{S}^{\lambda_N}a^N].
\end{align*}
\end{definition}

Exemplary, for $M=N=3$ and $\lambda=(1,-1,2)$ we obtain
\begin{align*}
S_{\lambda}\begin{pmatrix}
1 & 1 & 1 \\
2 & 2 & 2 \\
3 & 3 & 3
\end{pmatrix}=
\begin{pmatrix}
3 & 2 & 2 \\
1 & 3 & 3 \\
2 & 1 & 1
\end{pmatrix}.
\end{align*}
We observe that different values of $\lambda_k$ may produce the same shift. We can state some basic properties of shift operators:

\begin{theorem}\label{thmShiftBasic}
For arbitrary $\lambda\in\Z^N$ the operator $S_{\lambda}$ are linear. The inverse operators are given by $S_{\lambda}^{-1}=S_{-\lambda}$. Furthermore
$\|S_{\lambda}A\|_F=\|A\|_F$ holds. Given $\lambda'\in\Z^N$ we obtain
$S_{\lambda}S_{\lambda'}=S_{\lambda+\lambda'}$.
\end{theorem}

Using the shift operator, we define shifted matrices:

\begin{definition}[Shifted rank-$1$ matrices] We call $A\in\mathbb{C}^{M\times N}$ a shifted rank-$1$ matrix, if $A=S_{\lambda}(uv^*)$ for some vectors $u\in\mathbb{C}^M$, $v\in\mathbb{C}^N$ and $\lambda\in\mathbb{Z}^N$.
\end{definition}

\begin{remark}\label{rmk2}
Some full rank matrices such as the identity or convolution matrices can be written as shifted rank-$1$ matrix. Indeed, as already indicated in the introduction, each shifted rank-$1$ matrix can be written as product of two matrices $S_{\lambda}(uv^*)=UV$. Here $U$ is a convolution matrix constructed using the vector $u$. $V$ has $1$-sparse columns with non-zero elements $v$ at row positions $\lambda$. Using this factorization multiplying by a shifted rank-$1$ matrix can be performed in $O(M\log M+N)$ handling the convolution in Fourier domain.
\end{remark}

In Section  \ref{Sec::Greedy} we need the representation of shift operators in Fourier domain. Therefore, we denote by $F(A)=\hat{A}$ the column-wise Fourier transform of $A$. As a direct consequence of the Fourier shift theorem we obtain
\begin{align*}
F(S_{\lambda}A)=\hat{A}\odot P_{\lambda} && \text{with} &&
P_{\lambda}=\left[e^{-\frac{2\pi i}{M}j\lambda_k}\right]_{j,k=1}^{M,N}.
\end{align*}
Throughout this work, $P$ is a phase matrix, i.e., all entries have absolute value 1. By $P_{\lambda}$ we denote the specific phase matrix defined above.

\section{Shifted rank-1 approximation}\label{Sec::Algorithm}

This section describes the algorithm to find an approximate solution of (\ref{multishift}). The heuristic idea and advantages for some applications are discussed.

First, we note that (\ref{multishift}) is hard to solve for all variables at once. Hence, we develop an iterative technique where the number of unknowns in each step is drastically reduced. In each iteration only one of the $L$ unknown shifted rank-$1$ matrices is recovered. Afterwards the data $A$ is updated by subtracting the shifted rank-$1$ matrix. This process can be iterated $L$ times to reconstruct all unknowns. If $L$ is not given, the algorithm may stop when the remaining data  drops below some threshold.

Let us have a closer look at the iteration. In each step we seek the solution of 
\begin{align}\label{probApprox}
\min\limits_{\lambda,u,v}\|A-S_{\lambda}(uv^*)\|_F.
\end{align}
Note that this problem has infinite ambiguities. Given $S_{\lambda}(uv^*)$ other solutions are
\begin{align*}
S_{\lambda-m\mathbb{1}}(\tilde{S}^m(u)v^*), && S_{\lambda+(0,\ldots,0,M,0,\ldots,0)}(uv^*) \\
S_{\lambda}((\frac{u}{c}(cv^*)) && \text{ with }c\neq0,m\in\mathbb{Z}.
\end{align*}
To overcome these ambiguities we set $v_1=1$, $\lambda_0=0$ and $\lambda\in\{0,\ldots,M-1\}^{N}$ during this work. However, more ambiguities can occur e.g., if columns of $A$ are periodic.
Next, we rewrite problem (\ref{probApprox}) by applying the inverse shift operator and using its linearity
\begin{align*}
\min\limits_{\lambda,u,v}\|S_{-\lambda}(A)-uv^*\|_F.
\end{align*}
It turns out that once we know the shift vector $\lambda$, $u$ and $v$ form the standard rank-$1$ approximation which can be found by SVD. Hence, it remains the problem of finding the optimal $\lambda$. Let $\sigma_1,\ldots,\sigma_K$ be the singular values of $S_{-\lambda}(A)$. Then
\begin{align*}
\|S_{-\lambda}(A)-uv^*\|_F^2&=\sum\limits_{k=2}^K\sigma_k^2=\|S_{-\lambda}(A)\|_F^2-\sigma_1^2\\
&=\|A\|_F^2-\|S_{-\lambda}(A)\|_2^2.
\end{align*}
Thus, we can find $\lambda$ by solving 
\begin{align}\label{discreteMaxProb}
\max\limits_{\lambda\in\Z^N}\|S_{-\lambda}(A)\|_2^2.
\end{align}

Unfortunately, maximizing the spectral norm over a discrete set is still a hard problem. Note that even for $\lambda_j\in\{0,\ldots,M-1\}$ the number of possible vectors is $M^N$.  Thus, a naive combinatorial approach would lead to an exponential runtime. However, in the next section we introduce our approach to find a good approximation of (\ref{discreteMaxProb}). This method is a bit more technical and is applied in Fourier domain. Hence let us for now assume that we are able to solve (\ref{discreteMaxProb}). Then the shifted rank-$1$ approximation algorithm can be summarized as follows.

\begin{framed}
	\begin{algorithm}
		\begin{algorithmic}[1]
			\STATE{Input: $A$, $L$}
			\FOR{$k=1,\ldots,L$}
			\STATE{Solve (\ref{discreteMaxProb}) for $\lambda^k$}
			\STATE{Calculate $u^k$, $v^k$ using SVD of $S_{-\lambda}(A)$}
			\STATE{Update $A\leftarrow A-S_{\lambda}(u^k(v^k)^*)$}
			\ENDFOR
			\RETURN{$\lambda^k$, $u^k$, $v^k$ for $k=1,\ldots,L$}
		\end{algorithmic}
		\label{alg:sr1}
	\end{algorithm}
\end{framed}

Having the algorithmic scheme at hand, let us engross the idea behind this method. Given $A$ the algorithm can extract $L$ unique features that can be described as shifted rank-$1$ matrices. The obtained parameters $\lambda^k$, $u^k$ and $v^k$ may directly be used to extract information about these features without the need of any post-processing method. Thus, whenever the shifted rank-$1$ model applies to an application, the method can directly give crucial information about underlying features. Exemplary, let us demonstrate this for three applications.

Seismic exploration seeks for subsurface deposits of e.g., gas or oil. Therefore, an artificial seismic wave is generated and multiple sensors are placed along the testing area. In this setup earth layers play an important role as their boundaries reflect the seismic wave. The reflected wave varies depending on the material of the layer. Moreover, the time of arrival of the reflection at different sensors depends on the distance to the layer. A decomposition of the seismic image as given by our algorithm can be extremely useful. The shift vectors $\lambda^k$ give information about the position of each layer while the seismic wave $u^k$ (and also $v^k$) contain information about the material.

The non-destructive testing method named Time-of-Flight-Diffraction uses ultrasonic signals to detect defects inside materials such as steel tubes. Again, separating the data into $L$ unique features will give information about the position ($\lambda^k$) and type ($u^k$,$v^k$) of each feature. Automated processing can directly use this information to decide whether a tube needs to be fixed or not.

As a last example, we consider video analysis. We can represent a video as matrix $A$ where each frame is encoded as a column. Assume we have given a record of a traffic camera. These videos contain moving objects such as cars, cyclists or pedestrians. Each of these objects appears once per frame and hence is a unique event. Moreover, it changes its position in each frame while it moves through the image. Thus, the shifted rank-$1$ model applies and we can use the algorithm to track the movement of each object.

\section{Reconstruction of $\lambda$}\label{Sec::Greedy}

In this section we present an approach to find an approximate solution of (\ref{discreteMaxProb}). As noted before, a naive approach to solve this problem leads to exponential runtime. Even restricting the permitted shifting distance will not bring any advantage as the problem is exponential in the number of columns $N$. Another heuristic justification not to restrict the shift distance is given by the rank inequality of the Hadamard product
\begin{align}\label{HadIneq}
\begin{aligned}
&\mathop{rank}(A)=\mathop{rank}(\hat A)=\mathop{rank}(\hat A\odot P_{-\lambda}\odot P_{\lambda}) \\
\leq&\mathop{rank}(\hat A\odot P_{-\lambda})\mathop{rank}(P_{\lambda})=\mathop{rank}(S_{-\lambda}A)\mathop{rank}(P_{\lambda}).
\end{aligned}
\end{align}
The rank of $A$ is typically large, otherwise a normal low-rank approximation could be applied. Hence, if we want to achieve the best case $\mathop{rank}(S_{-\lambda}A)=1$, we need  $\mathop{rank}(P_{\lambda})$ to be large. From the definition of $P_{\lambda}$ it directly follows that the rank of this matrix coincides with the number of different shifts in $\lambda$.

Inequality (\ref{HadIneq}) already demonstrates that stating the problem in Fourier domain may give some insights that can not be achieved otherwise. In most parts of this section we use the representation in Fourier domain. This brings two advantages. First, problem (\ref{discreteMaxProb}) can be relaxed in Fourier domain. The relaxed problem has no adequate formulation in time domain but its maxima can be described easily. Second, performing the algorithm in Fourier domain decreases the runtime to $O(NM\log M)$ since it basically evaluates cross-correlations that transfer to element-wise multiplications in Fourier domain.

Let us start with a simple method that can find a local maximum of (\ref{discreteMaxProb}). Therefore let $u$ be the first left singular vector of $A$ and $a^1,\ldots,a^N$ be the columns of $A$, i.e., $A=[a^1,\ldots,a^N]$. Suppose we have given a shift vector $\lambda$ such that $\|u^*A\|_2<\|u^*S_{-\lambda}(A)\|_2$ holds. Then it follows that
\begin{align}\label{incLambda}
\|A\|_2=\|u^*A\|_2<\|u^*S_{-\lambda}(A)\|_2\leq\|S_{-\lambda}A\|_2.
\end{align}
Hence, the spectral norm increases. To ensure the strict inequality holds note that
\begin{align}\label{equalDiff}
\|u^*S_{-\lambda}(A)\|_2^2-\|u^*A\|_2^2=\sum\limits_{k=1}^N\left(\left|u^*\tilde{S}^{-\lambda_k}a^k\right|^2-|u^*a^k|^2\right).
\end{align}
Note that in this sum all coefficients of $\lambda$ can be optimized independently to achieve the maximum difference between both values. However, in our numerical experiments it turned out that updating all entries of $\lambda$ at once leads to an unstable method that often runs into small local maxima. Hence, we solve the following problem.
\begin{align}\label{localMax}
\max\limits_{k, \lambda^k}\left|u^*\tilde{S}^{-\lambda_k}a^k\right|^2-|u^*a^k|^2.
\end{align}
There are two possible scenarios. First, the maximum might be zero. It follows that no choise of $\lambda_k$ for any column $k$ can increase the spectral norm and thus we converged to a local maximum. Second, the maximum is larger than zero. Let $k$ and $\lambda_k$ be the position where the maximum is achieved. Then following from (\ref{incLambda},\ref{equalDiff}) shifting the $k$-th column by $-\lambda_k$ will increase the spectral norm. We update $A$ accordingly, calculate the new first left singular vector $u$ and repeat this procedure until ending in a local maximum. This is guaranteed to happen after a finite number of steps since there are only finitely many possibilities for $\lambda$.

Note that the first term in (\ref{localMax}) asks for the maximum of the cross-correlation of $u$ and $a^k$. The cross-correlation can efficiently be calculated in Fourier domain as $F^{-1}(\overline{\hat{u}}\odot\hat{a^k})$. The second term in (\ref{localMax}) coincides with the first coefficient of the cross-correlation. We slightly abuse the notation by defining
$$\overline{\hat{u}}\odot\hat{A}:=[\overline{\hat{u}}\odot\hat{a^1},\ldots,\overline{\hat{u}}\odot\hat{a^N}]$$
as the element-wise product of a vector and each column of the matrix. Define $B=|F^{-1}(\overline{\hat{u}}\odot\hat{A})|^2$ and let $b^1$ be its first row. Then (\ref{localMax}) can be written as finding the maximum of $B-\mathbb{1}b^1$. We can summarize the local maximization as Algorithm \ref{alg:localOpt}.

\begin{framed}
	\begin{algorithm}
		\begin{algorithmic}[1]
			\STATE{Input: $A$}
			\STATE{Set $s=1$}
			\WHILE{$s\neq0$}
			\STATE{Calculate $\hat{u}$ the first left singular vector of $\hat{A}$}
			\STATE{Set $B=|F^{-1}(\overline{\hat{u}}\odot\hat{A})|^2$}
			\STATE{Subtract the first row of $B$ from every row}
			\STATE{Get the position $(k,s)$ of the maximum in $B$}
			\STATE{Shift the $k$-th column of $A$ by $-s$}
			\STATE{Update $\lambda_k\leftarrow\lambda_k+s$}
			\ENDWHILE
			\RETURN{$\lambda$}
		\end{algorithmic}
		\label{alg:localOpt}
	\end{algorithm}
\end{framed}

\begin{remark}\label{rmk1}
Algorithm \ref{alg:localOpt} requires the computation of element-wise multiplications,  (inverse) Fourier transforms and the calculation of the first left singular vector of $\hat{A}$. The multiplications scale linear in $M$ and $N$ and the inverse Fourier transform can be calculated in $O(MN\log M)$. Computing a full singular value decomposition of $\hat{A}$ would have complexity $O(\min(M^2N,MN^2))$. However, since we only need the first left singular vector, much more efficient algorithms can be used. One possibility is, to use effective methods for rank-$1$ updates of the matrix, since shifting one column is a rank-$1$ update. The computational costs than reduce to $O(MN)$ \cite{updateSVD}. As another approach, we can use the previous singular vector as a starting guess for any iterative method, e.g., power iteration. Since the matrix was only shifted in one column, this will be a very good guess. Altogether, the local optimization has a complexity of $O(MN\log M)$ per iteration. The same arguments apply for the following techniques.

Moreover, since most of the method is performed in Fourier domain. The Fourier transforms only need to be calculated once. Shifting the matrix $A$ can also be implemented in Fourier domain by multiplying with a phase shift.
\end{remark}

Using Algorithm \ref{alg:localOpt} we are able to find a local maximum of (\ref{discreteMaxProb}). Because we can only guarantee local convergence, the choice of the starting guess is crucial. Especially in the presence of noise a wrong starting point may end up in a small local maximum. Note that problem (\ref{discreteMaxProb}) is indeed stable to Gaussian noise. Let $E$ be random Gaussian noise. Using that Gaussian noise is shift-invariant we obtain 
	\begin{align*}
	\|S_{-\lambda}(A+E)\|_2=\|S_{-\lambda}(A)+E\|_2\approx\|S_{-\lambda}A\|_2,
	\end{align*}
where the last approximation holds because Gaussian noise corrupts all singular values equally. In the remainder of this section we introduce a strategy that, according to our numerical experiments, leads stable results independent of the type of input data. Therefore, we formulate problem (\ref{discreteMaxProb}) in Fourier domain as
\begin{align}\label{discreteMaxProbFourier}
\max\limits_{\lambda\in\Z^N}\|\hat{A}\odot P_{-\lambda}\|_2^2.
\end{align}
Now consider a relaxed version of (\ref{discreteMaxProbFourier}) where the dependency on $\lambda$ is neglected and any phase matrix can be used
\begin{align}\label{relaxedFourier}
\max\limits_{P,|p_{jk}|=1}\|\hat{A}\odot P\|_2^2.
\end{align}
We state the following theorem.
\begin{theorem}\label{relaxSolution}
It holds that 
\begin{align*}
\||\hat{A}|\|_2^2=\max\limits_{P,|p_{jk}|=1}\|\hat{A}\odot P\|_2^2.
\end{align*}
Let $\hat{u}^\text{opt}$ be the first left singular vector of $|\hat{A}|$ and $\hat{u}^P$ the first left singular vector of $\hat{A}\odot P$. Then
\begin{align*}
\|\hat{A}\odot P\|_2=\||\hat{A}|\|_2 && \Rightarrow && |\hat{u}^P|=|\hat{u}^\text{opt}|.
\end{align*}
Furthermore, if $A=S_{\lambda}(uv^*)$ is a shifted rank-$1$ matrix, then $|\hat{u}^\text{opt}|=|\hat{u}|=|\hat{u}^{P_{-\lambda}}|$.
\end{theorem}
\begin{IEEEproof}
The first statement follows from
\begin{align*}
\|\hat{A}\odot P\|_2^2=\|(\hat{u}^P)^*(\hat{A}\odot P)\|_2^2
\leq\||\hat{u}^P|^*|\hat{A}|\|_2^2
\leq\||\hat{A}|\|_2^2.
\end{align*}
Note that the second inequality only holds if $|\hat{u}^P|$ is the first left singular vector of $|\hat{A}|$. Using the Perron-Frobenius theorem this is equivalent to $|\hat{u}^P|=|\hat{u}^\text{opt}|$.

For $A=S_{\lambda}(uv^*)$ we have $|\hat{A}|=|\hat{u}v^*|=|\hat{u}||v^*|$ and thus $|\hat{u}^\text{opt}|=|\hat{u}|$. The second equality follow from $\hat{A}\odot P_{-\lambda}=\hat{u}v^*$.
\end{IEEEproof}
Theorem \ref{relaxSolution} gives us the maximum value of the relaxed problem (\ref{relaxedFourier}) and states a necessary condition. Furthermore, it shows that this necessary condition can be fulfilled for shifted rank-$1$ matrices using the correct shift. We use these information to estimate a starting guess for our local maximization algorithm. The starting guess is calculated in two steps. First, we try to get close to the optimal value $\||\hat{A}|\|_2^2$. Second, we iterate towards the necessary condition. Both the singular vectors and values continuously depend on the matrix coefficients. Taking a starting guess close to the necessary condition of the relaxed problem and with a large first singular value, we may be close to a global optimum of (\ref{relaxedFourier}) and hence converge to a large local maximum of (\ref{discreteMaxProbFourier}).

Let $\Omega\in\C^{N\times N}$ be any diagonal matrix with $|\omega_{kk}|=1$. Then $\||\hat{A}|\|_2=\||\hat{A}|\Omega\|_2$. We construct our starting guess by pulling the largest singular value of $\hat{A}\odot P_{-\lambda}$ towards this value. Therefore consider the inequality

\begin{align}\label{startguessError}
\begin{aligned}
&\||\hat{A}|\Omega\|_2^2-\|\hat{A}\odot P_{-\lambda}\|_2^2\\
\leq&\||\hat{A}|\Omega-\hat{A}\odot P_{-\lambda}\|_2^2\\
\leq& \||\hat{A}|\Omega-\hat{A}\odot P_{-\lambda}\|_F^2\\
=&2(\|A\|_F^2-\mathop{Re}\langle|\hat{A}|\Omega,\hat{A}\odot P_{-\lambda}\rangle)
\end{aligned}
\end{align}
where $\mathop{Re}\langle\cdot,\cdot\rangle$ is the real part of the inner product. The columns of the involved matrices are given as
\begin{align*}
|\hat{A}|&=[|\hat{a}^1|,\ldots,|\hat{a}^N|],\\
\hat{A}\odot P_{-\lambda}&=[F(\tilde{S}^{\lambda_1}a^1),\ldots,F(\tilde{S}^{\lambda_N}a^N)].
\end{align*}
The matrix inner product in (\ref{startguessError}) can be formulated as sum of inner products of the columns. We obtain
\begin{align*}
\mathop{Re}\langle|\hat{A}|\Omega,\hat{A}\odot P_{-\lambda}\rangle
&=\mathop{Re}\sum\limits_{k=1}^{N}\overline{\omega}_{kk}\langle|\hat{a}^k|,F(\tilde{S}^{\lambda_k}a^k)\rangle\\
=\sum\limits_{k=1}^{N}\left|\langle |\hat{a}^k|,F(\tilde{S}^{\lambda_k}a^k)\rangle\right|
&=\sum\limits_{k=1}^{N}\left|\langle F^{-1}(|\hat{a}^k|),\tilde{S}^{\lambda_k}a^k\rangle\right|
\end{align*}
where we choose $\omega_{kk}$ to be the phase of the inner product. Minimizing (\ref{startguessError}) can be done by maximizing each inner products. Similar as in (\ref{localMax}), maximizing over $\lambda$ coincides with finding the maximum of the cross-correlation of the involved vectors $F^{-1}(|\hat{a}^k|)$ and $a^k$. Hence we use the representation of the cross-correlation in Fourier domain an calculate the matrix $B=|F^{-1}(|\hat{A}|\odot\hat{A})|^2$. The starting guess $\lambda$ is given by the row index of the maximum value of each column of $B$.

Before we iterate to a local maximum, we use a modified version of Algorithm \ref{alg:localOpt} to force the largest singular vector towards the necessary condition stated in Theorem \ref{relaxSolution}. After the singular vector has been calculated in Line 3 of the algorithm, we add an additional step to integrate the optimal amplitude, i.e., we update
\begin{align}\label{modAlg2}
\hat{u}\leftarrow|\hat{u}^\text{opt}|\odot\mathop{phase}(\hat{u}).
\end{align}
Note that using this modification it is not guaranteed that Algorithm \ref{alg:localOpt} will terminate in a local maximum of (\ref{discreteMaxProb}). Hence, local optimization is applied afterwards. Altogether we obtain the following algorithm.

\begin{framed}
	\begin{algorithm}
		\begin{algorithmic}[1]
			\STATE{Input: $A$}
			\STATE{Calculate $B=|F^{-1}(|\hat{A}|\odot\hat{A})|^2$}
			\STATE{Set $\lambda$ to the indices of the maximum values of $B$ in each column and update $A\leftarrow S_{-\lambda}A$}
			\STATE{Update $\lambda$ and $A$ using Algorithm \ref{alg:localOpt} with (\ref{modAlg2})}
			\STATE{Update $\lambda$ and $A$ using Algorithm \ref{alg:localOpt}}
			\RETURN{$\lambda$}
		\end{algorithmic}
		\label{alg:slr}
	\end{algorithm}
\end{framed}
This subroutine is used in Algorithm \ref{alg:sr1} at Line 3 and thus completes our approximation method. Following Remark \ref{rmk1} the runtime of Algorithm \ref{alg:sr1} is $O(LMN\log M)$.

\section{Numerical results}

\begin{figure}
\centering
\includegraphics[width=0.24\textwidth]{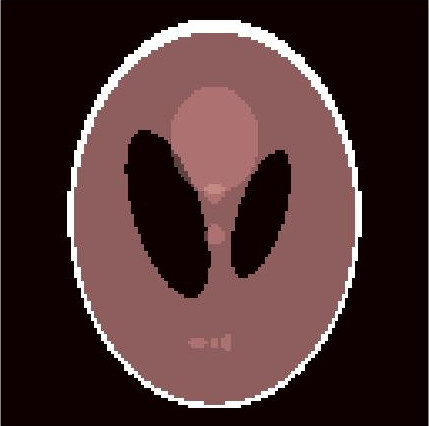}
\includegraphics[width=0.24\textwidth]{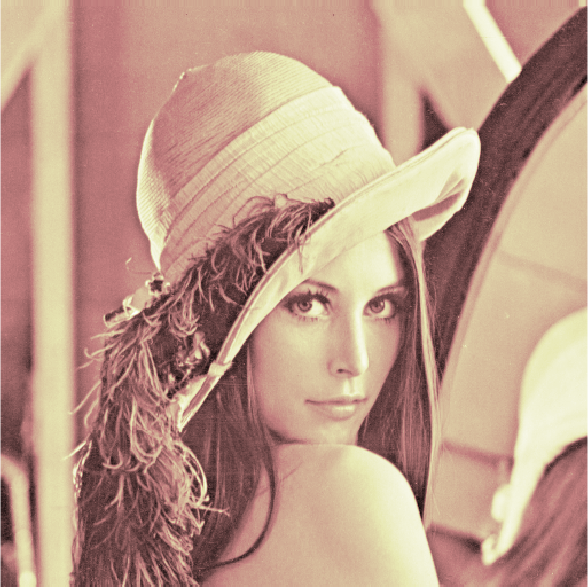}\\[0.1cm]
\includegraphics[width=0.24\textwidth]{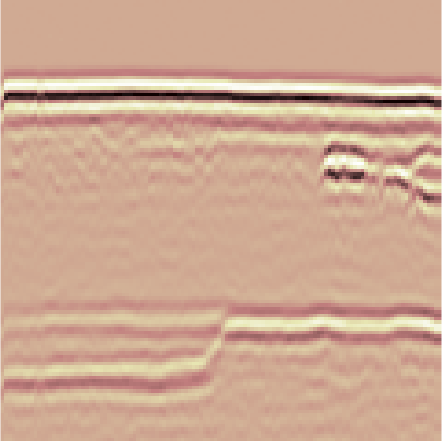}
\includegraphics[width=0.24\textwidth]{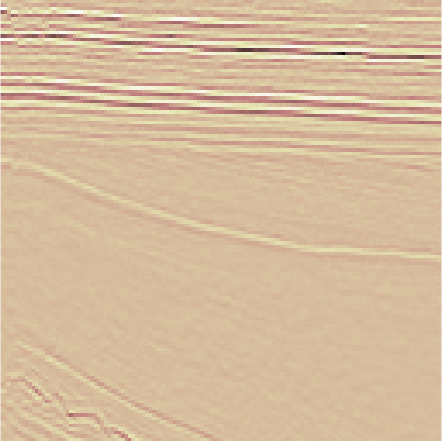}
\caption{Exemplary input data: Cartoon-like, natural, ultrasonic and seismic images.}
\label{fig:ExampleData}
\end{figure}

In this section we analyze our methods in several numerical experiments divided in five subsections. In the first subsection we investigate the general performance of the method and the shift approximation strategy. The second subsection compares the approximation quality to four other methods: SVD, wavelet transform, MoTIF and UC-DLA. The latter two are of special interest as they are close to the technique discussed in this work. The runtime of our technique is discussed in a third subsection. In all three subsections we apply the methods to different kinds of data with $M=N=128$, some exemplary shown in Figure \ref{fig:ExampleData}. We use cartoon-like, natural, seismic, and ultrasonic images as well as random orthogonal matrices. The Frobenius norm of all input data is normalized to $1$. Cartoon-like and natural images are chosen as standard test for image processing algorithms. Seismic and ultrasonic images fit into the shifted rank-$1$ model. Random orthogonal matrices  are usually hard to approximate. Note that the illustrated figures are created using real data. However, in the experiments we also applied the method to complex data obtaining similar results.

The last two subsections demonstrate the benefits of our method in applications. We use the algorithm for segmentation of highly noised data. Furthermore, we track objects in videos such as soccer match recordings.

\subsection{Performance behavior}

\begin{figure}
\centering
\subfloat[]{\includegraphics[width=0.24\textwidth]{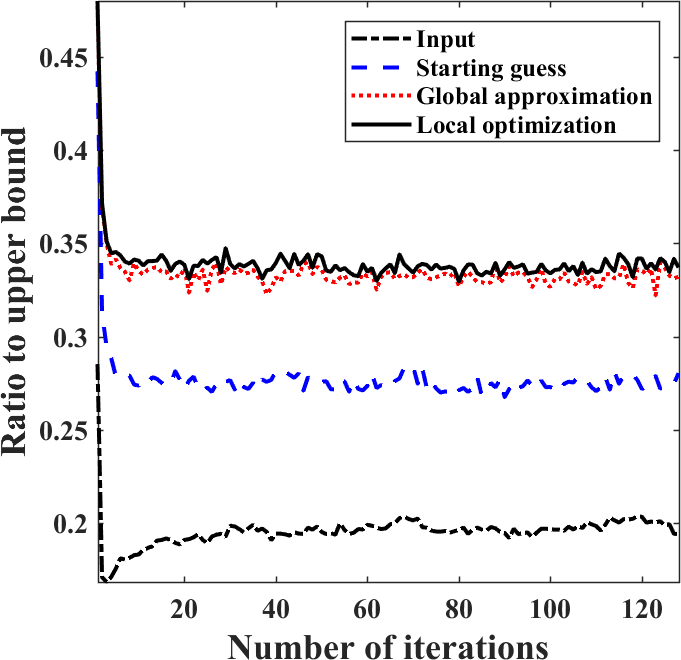}\label{fig:sv_ratio}}
\subfloat[]{\includegraphics[width=0.24\textwidth]{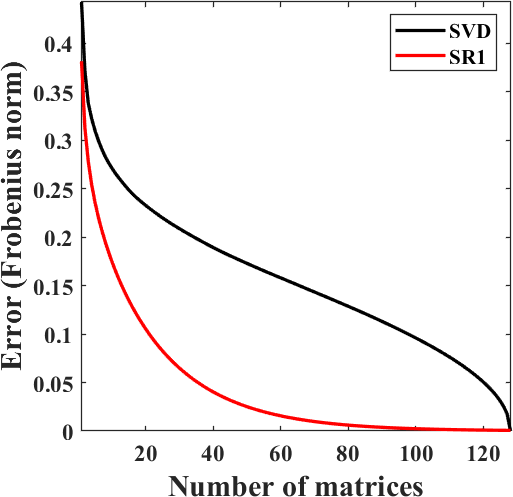}
\label{fig:ErrUnited}}
\caption{(a) Singular value ratio to $\||\hat{A}|\|_2$ after the individual steps of Algorithm \ref{alg:slr}.
(b) Average approximation error over all kinds of input data.}
\end{figure}

In this subsection we give a first demonstration of the general approximation performance of our proposed method which we will shortly name SR1 (shifted rank-$1$). In Figure \ref{fig:ErrUnited} we plotted the average approximation error over all data used in the experiments of SR1 against the number $L$ of matrices used. As comparison, we included the approximation error of low-rank approximation using SVD. We observe that, in contrast to SVD, SR1 seems to have exponential decay independent of the underlying data.

Algorithm \ref{alg:slr} derives the shift $\lambda$ in three steps. First, a starting guess is calculated (Line 3). Afterwards, it tries to approximate a global solution in Line 4. Last, local optimization is applied. We clarify the importance of every step in Figure \ref{fig:sv_ratio}. Here, the ratio of the largest singular value to the upper bound $\||\hat{A}|\|_2$ for the input data and after each of the three steps is plotted. Clearly, starting guess and global approximation greatly improve the results. Moreover, note that the obtained ratio does not change with the number of iterations $L$, i.e., the algorithm does not suffer from saturation. Local optimization only causes a minor improvement. Table \ref{tab:iterations} gives an explanation for this. Here, the number of iterations performed during Algorithm \ref{alg:localOpt} is shown. The calculated shift is already quite optimal such that local optimization terminates after only a few iterations. Remember that in each iteration only one column is shifted and the tested data is of size $M=N=128$. Hence for cartoon-like images, seismic or ultrasonic data the mean number of all shifts is even less than the number of columns.

\begin{table}
	\centering
	\caption{Mean number of iterations for different kind of input data.}
	\begin{tabular}{|l|c|c|c|}
		\hline
		data & global & local & total\\
		\hline
		orthogonal & 139 & 21 & 160\\
		natural & 137 & 24 & 161 \\
		cartoon & 91 & 15 & 106 \\
		seismic & 95 & 16 & 111 \\
		ultrasound & 95 & 18 & 113 \\
		\hline
	\end{tabular}
	\label{tab:iterations}
\end{table}

\begin{figure}
\centering
\includegraphics[width=0.15\textwidth]{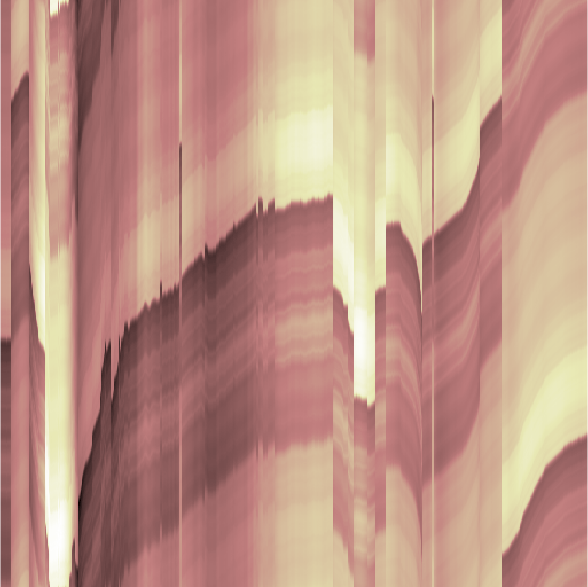}
\includegraphics[width=0.15\textwidth]{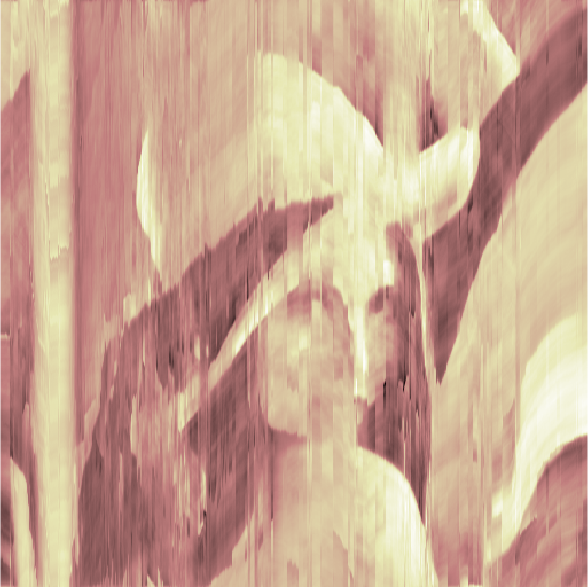}
\includegraphics[width=0.15\textwidth]{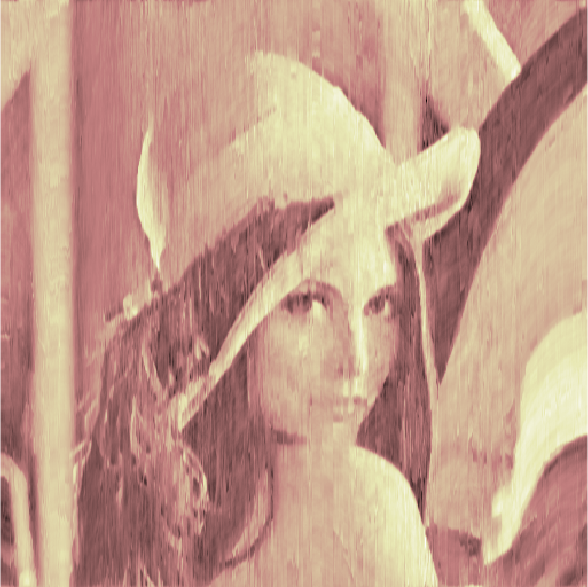}
\caption{Reconstruction of Lena image using 1, 5 and 10 shifted rank-$1$ matrices.}
\label{fig:Lena}
\end{figure}

Last in this subsection, we want to illustrate the influence of the shift vector on the approximation. Therefore we computed the shifted rank-$1$ approximatino of the Lena image using $1$, $5$ and $10$ matrices. The result is shown in Figure \ref{fig:Lena}. As we see, by shifting the columns we are able to follow horizontal structures in the image. This becomes very clear looking at the most left image where only one shifted rank-$1$ matrix is used. Vertical structures such as the hair cannot be reconstructed in such detail using column shifts (see middle image).

\subsection{Sparse approximation}

\begin{figure}
\centering
\includegraphics[width=0.24\textwidth]{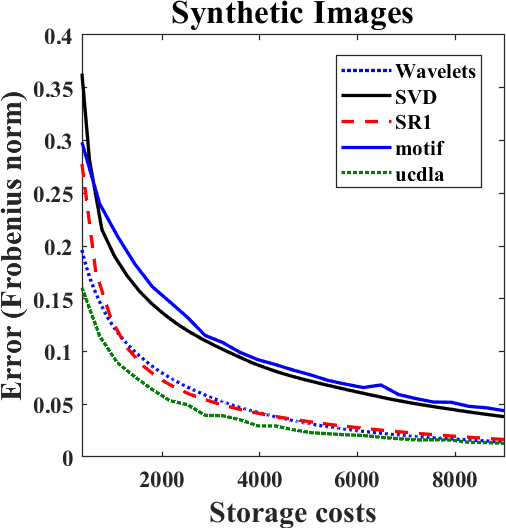}
\includegraphics[width=0.24\textwidth]{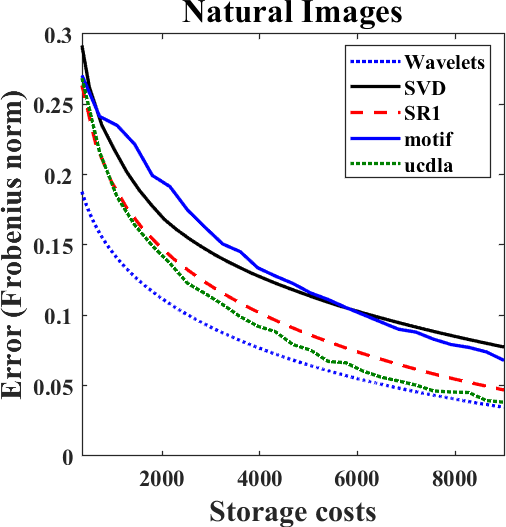}
\includegraphics[width=0.24\textwidth]{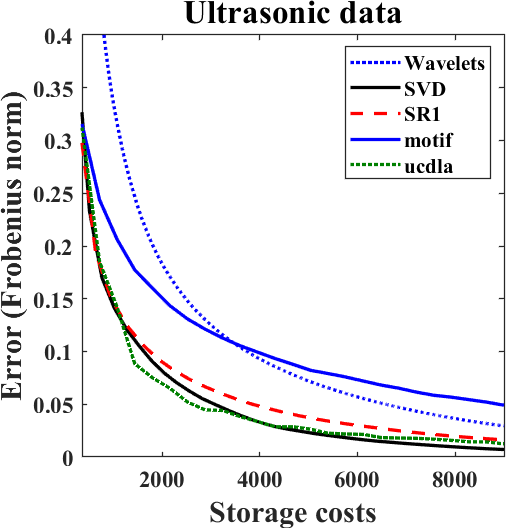}
\includegraphics[width=0.24\textwidth]{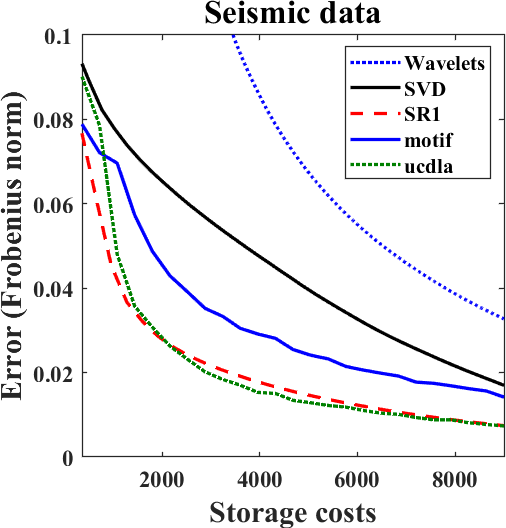}
\includegraphics[width=0.24\textwidth]{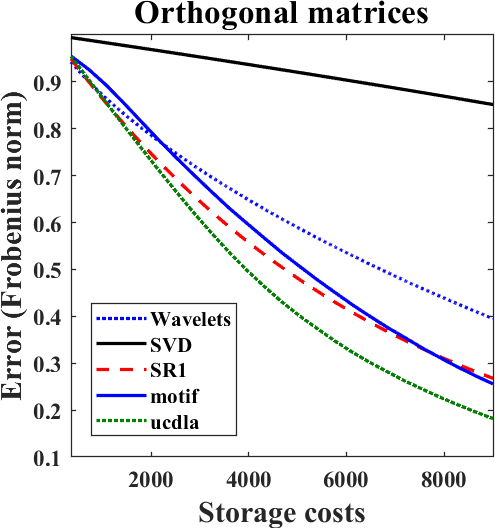}
\caption{Approximation error of all algorithms for different kinds of input data plotted against the storage costs.}
\label{fig:storage}
\end{figure}

Now, we analyze the storage costs of our method compared to the obtained approximation accuracy. We compare the results with four methods: low rank approximation using SVD , Daubechies7 wavelets, MoTIF and UC-DLA. Storing one shifted rank-$1$ matrix requires $M+N$ doubles and $N$ integers which we count as $0.5$ doubles. We identify the storage costs of a double scalar as $1$. Hence, in this experiment with $M=N=128$ we have a storage costs of $320$ per rank-$1$ matrix. Storage costs for the other methods are counted similarly. The obtained approximation error is shown in Figure \ref{fig:storage}.

We note that SR1 is for both ultrasonic and seismic images one of the best performing methods. For these type of images the shifted rank-$1$ model fully applies. SR1 achieves similar approximation results as UC-DLA even though it has a more restrictive model. MoTIF approximation accuracy suffers from its uncorrelated dictionary design. Typical atoms in ultrasonic or seismic data will be highly correlated. Moreover, the underlying testing setup of the ultrasonic images promotes low-rank data. Hence SVD performs especially good.

\begin{figure}
	\centering
	\includegraphics[width=0.15\textwidth]{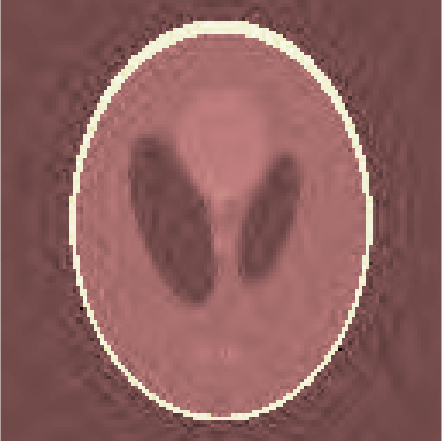}
	\includegraphics[width=0.15\textwidth]{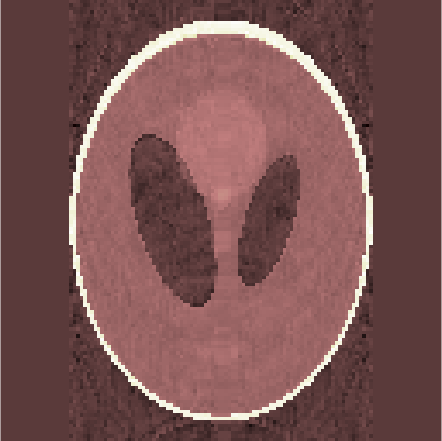}
	\includegraphics[width=0.15\textwidth]{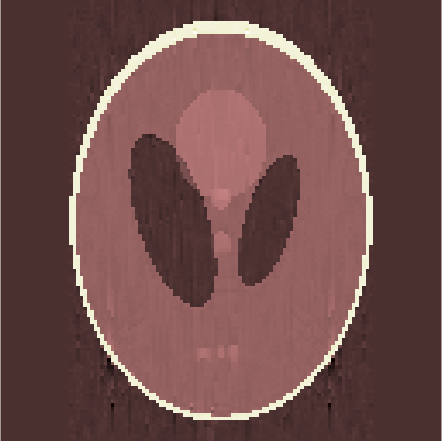}
	\caption{Sparse approximation of image "phantom" using Wavelets (left), SR1 (middle) and UC-DLA (right).}
	\label{fig:synEx}
\end{figure}

Still, wavelet approximation is the best of the chosen methods for natural images which can contain many details for which the shift approach is not optimal. However, already for cartoon-like images, shifted rank-$1$ approximation is as good as wavelet compression. This is due to the fact that cartoon-like images contain more clear edges and smooth surfaces which can be shifted into a low-rank form. Consequently, UC-DLA outperforms SR1 and Wavelet transform slightly. The adapted shift-invariant dictionary can show its full strength for this type of data. Exemplary, Figure \ref{fig:synEx} shows the reconstructed image of Wavelets, SR1 and UC-DLA with a storage costs of $3200$ (i.e., $L=10$). Note that, although the approximation error of all three methods is similar, SR1 and UC-DLA perform much better in preserving edges of the image. We have already seen this phenomena in the last subsection for the Lena image.

All methods that include a shift-invariance perform good in approximating orthogonal matrices. Altogether, we can summarize that SR1 calculates a reasonable sparse approximation. Whenever the data fits into the shifted rank-$1$ assumption, the method is one of the best. Due to the restricting model the approximation quality might suffer slightly. However, this will pay of in application as later examples demonstrate.

\subsection{Runtime}

Im Remark \ref{rmk1} we stated that SR1 has a runtime of $O(LMN\log M)$. Furthermore, according to Remark \ref{rmk2} multiplying with a shifted rank-$1$ matrix can be significantly faster. We test both statements in numerical experiments shown in this subsection.

\begin{figure}
	\centering
	\subfloat[]{\includegraphics[width=0.24\textwidth]{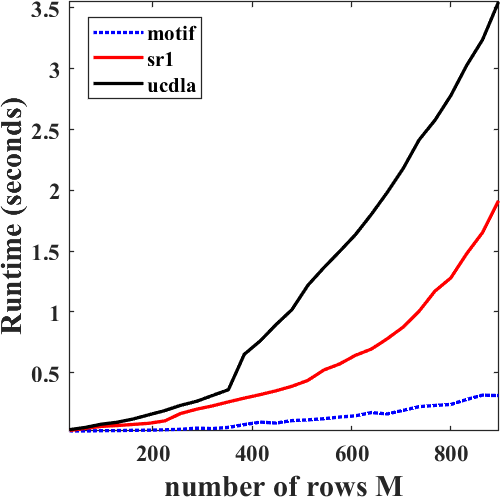}\label{fig:runtime}}
	\subfloat[]{\includegraphics[width=0.24\textwidth]{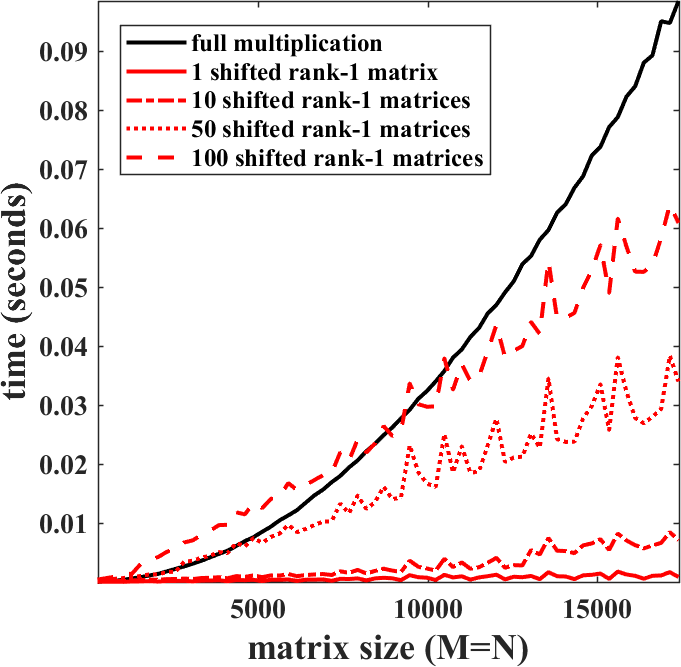}
		\label{fig:runtime_multiply}}
	\caption{(a) Average runtime of data approximation against the number of rows.
		(b) Average runtime of matrix vector multiplication using different number of shifted rank-$1$ matrices.}
\end{figure}

Let us first test the runtime of SR1. The linear scaling in $L$ is obvious since this is the number of iterations. Also a linear dependence in the number of columns $M$ is easy to see. Most operations of the algorithm are applied column-wise and thus could even run highly parallel. Most interesting is the runtime according to varying $M$. The optimal bound $O(M\log M)$ can only be achieved when the update of the singular vector is implemented in a clever way. For our implementation we use a power method which might not lead to optimal rates. We calculated the mean runtime for $N=128$ with increasing number $M$. We compare our method with the runtime of MoTIF and UC-DLA. The results are shown in Figure \ref{fig:runtime}. Although the rate seems worse than the theoretical optimum, we can outperform UC-DLA. Only the MoTIF algorithm, which showed worse approximation results in the last subsection, is much faster.

Second, we numerically investigate the fast matrix vector multiplication $Ax=b$ where $A$ is the sum of $L$ shifted rank-$1$ matrices. We measured the mean runtime for $L=1,10,50,100$ and different matrix sizes $M=N=256,512,\ldots,17408$. The results are shown in Figure \ref{fig:runtime_multiply}. As theoretically proven, the runtime is nearly linear and outperforms the full matrix multiplication even for high $L$. Note that the spikes that occur at certain matrix sizes coincide with values of $M$ that have a large prime factor and thus the fast Fourier transform is less efficient. Hence, zero padding to increase $M$ up to a power of $2$ may be a good option.

Both speed tests were performed in Matlab 2017a on Windows 10 with an i7-7700K CPU (4.2Ghz) and 32GB RAM.

\subsection{Segmentation}

\begin{figure}
	\centering
	\includegraphics[width=0.15\textwidth]{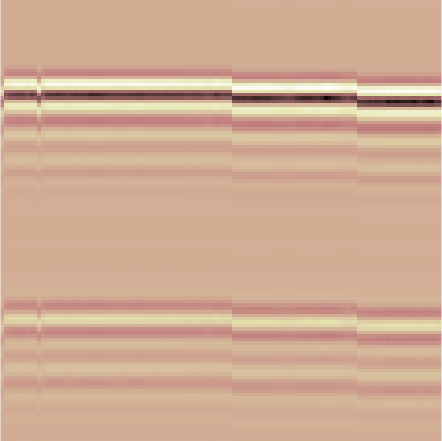}
	\includegraphics[width=0.15\textwidth]{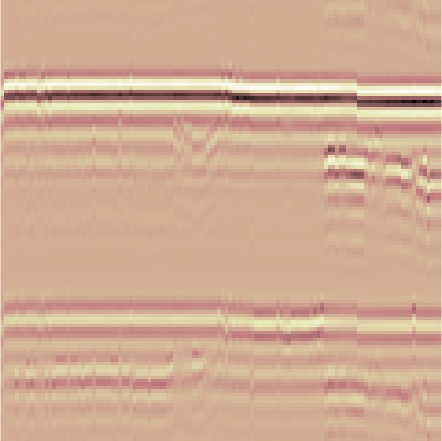}
	\includegraphics[width=0.15\textwidth]{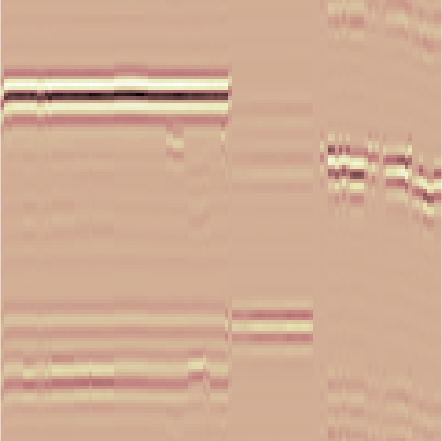}\\[0.1cm]
	\includegraphics[width=0.15\textwidth]{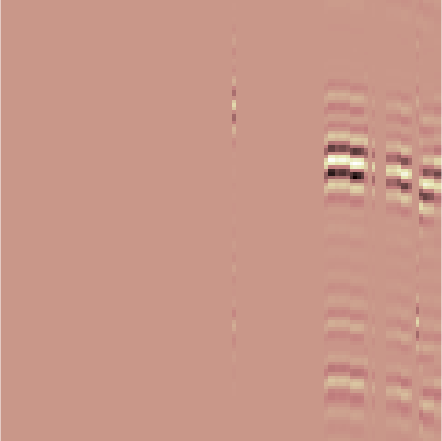}
	\includegraphics[width=0.15\textwidth]{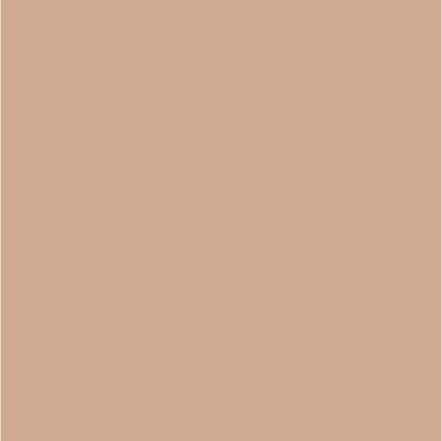}
	\includegraphics[width=0.15\textwidth]{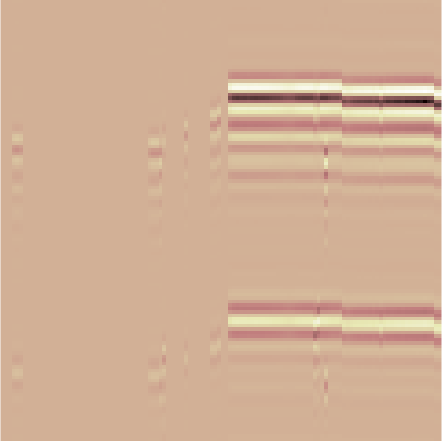}
	\caption{Separation of an ultrasonic image in two signals (top and bottom) using SR1 (left), MoTIF (middle) and UC-DLA (right).}
	\label{fig:ultraSep}
\end{figure}

Our first applicational demonstration comes from the field of non-destructive testing. The ultrasonic image shown in Figure \ref{fig:ExampleData} (bottom, left) was obtained measuring a defect in a steel tube. On the top and bottom of the image one can see the reflections of surface and back-wall of the tube. Another signal can be found at the right side of the image directly below the surface reflection. This signal indicates a defect in the material. We apply our method with $L=2$ to separate the defect from normal surface and back-wall reflections. We also use MoTIF and UC-DLA with a shift-invariant dictionary built from two atoms to separate the data. The results are shown in Figure \ref{fig:ultraSep}. The top row shows the first separated signal while the bottom row shows the second. MoTIF (middle) assigns both signals to the same atom. Also UC-DLA (right) fails to separate the signals. Only SR1 (left) is able to isolate the defect.

\begin{figure}
	\centering
	\includegraphics[width=0.15\textwidth]{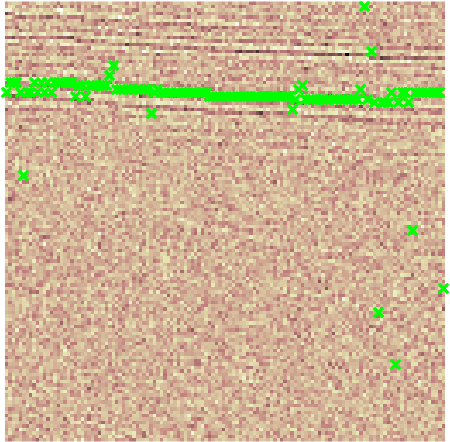}
	\includegraphics[width=0.15\textwidth]{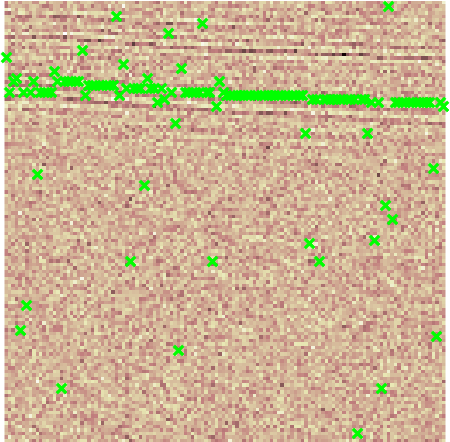}
	\includegraphics[width=0.15\textwidth]{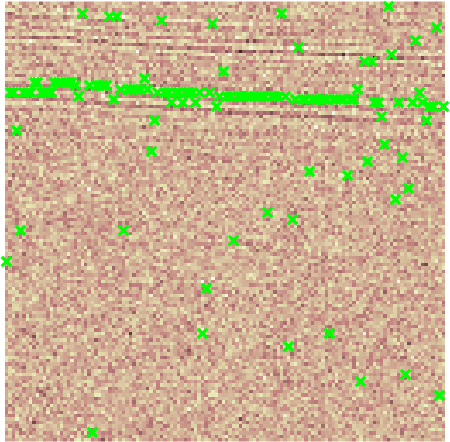}
	\caption{Identified earth layer reflection in noisy seismic image.}
	\label{fig:denoise}
\end{figure}

In a second test, we corrupt the seismic image shown in Figure \ref{fig:ExampleData} (bottom right) with noise. We will use our algorithm with $L=1$ to identify the first earth layer reflection at the top of the image. Figure \ref{fig:denoise} shows three images corrupted by Gaussian noise. The PSNR of the noisy images are $20$, $19$ and $18$ (left to right). The shift vector $\lambda^1$ is plotted on top of the images to indicate the identified earth layer. Even for high noise levels the reconstruction is still good.

\subsection{Object tracking}

In our last example we track objects in video records. The videos used here are contained in the "UCF Sports Action Data Set" downloaded from the UCF Center for Research in Computer Vision \cite{Video1,Video2}. Before we can apply our algorithm, we have to transform the 3D video data into a 2D matrix. Therefore, each frame will be transformed into a column vector.

\begin{figure}
	\centering
	\includegraphics[width=0.4\textwidth]{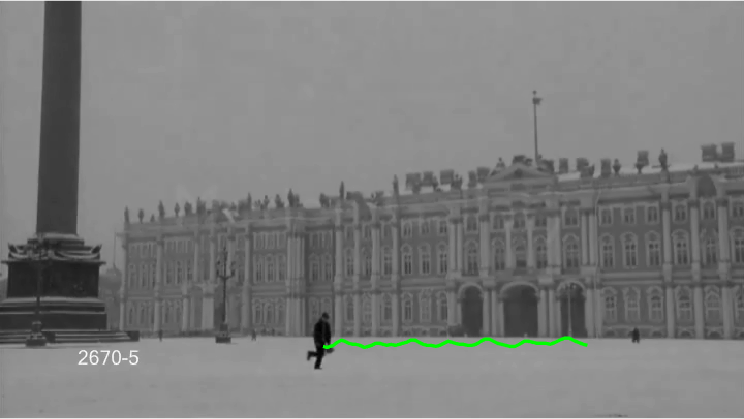}\\[0.1cm]
	\includegraphics[width=0.4\textwidth]{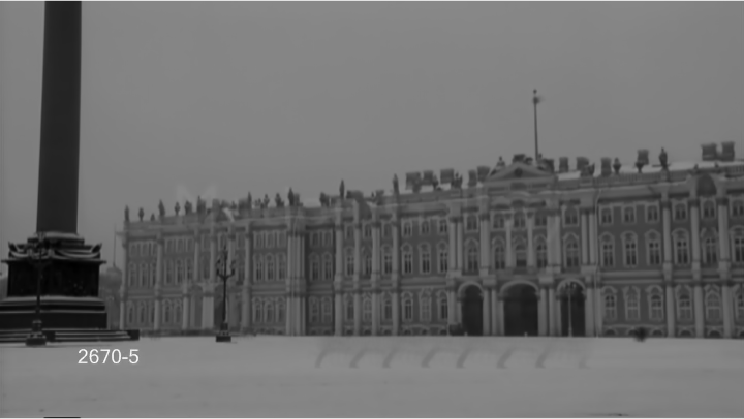}\\[0.1cm]
	\includegraphics[width=0.4\textwidth]{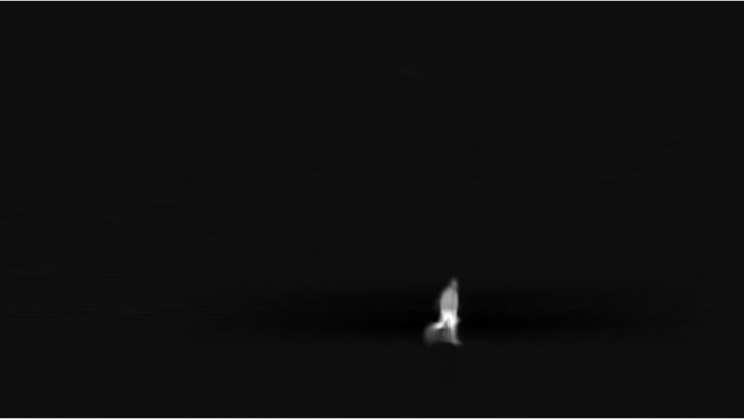}
	\caption{Tracked route in original video (top), reconstructed singular vectors $u^1$ (middle) and $u^2$ (bottom).}
	\label{fig:video1}
\end{figure}

We start with a simple example. The first video shows a person running from left to right with a constant background and constant camera position. We can extract the person by calculating the shifted rank-$1$ approximation with $L=2$. Now $u^k$ represents an object in one frame and $\lambda^k$ describes its movement. In Figure \ref{fig:video1} we have plotted both objects $u^1$ and $u^2$ (middle and bottom). As one can see, the constant background separates from the moving person. Some minor artifacts appear along the path of the person. This happens since our model can capture general movement (the persons path) but not the changing shape (moving arms and legs). The top image shows the first frame of the original video where we plotted the shift $\lambda^2$ to indicate the tracked route of the person.

\begin{figure}
	\centering
	\includegraphics[width=0.4\textwidth]{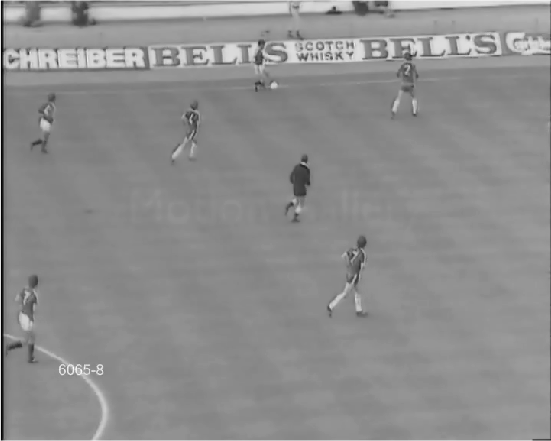}\\[0.1cm]
	\includegraphics[width=0.4\textwidth]{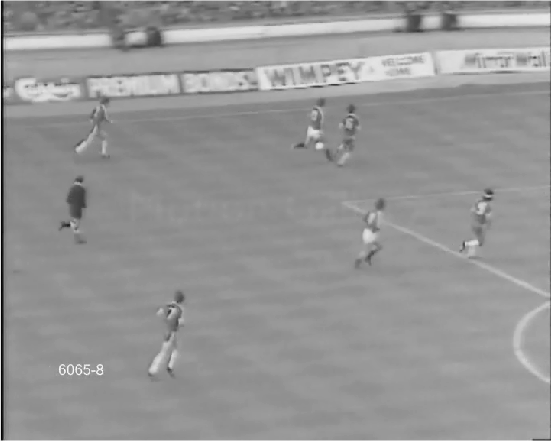}
	\caption{First and last frame of the soccer video clip.}
	\label{fig:video2}
\end{figure}

The next video is much more challenging. This clip shows a soccer game. Here, we have many moving objects such as players, referee or ball. Furthermore, the camera itself moves. This leads to a moving background where e.g., the advertisements at the sideline or the grass marks change. Figure \ref{fig:video2} shows the first and last frame of the video.

There are two main problems that do not fit in our setting. First, some players and part of the background move outside / run into the video what does not fit the circulant shifts. Second, the objects are not additive, i.e., their gray scale values do not add up but always the value of the front object is shown. This can lead to artifacts especially when the objects have similar gray scale values such as players and field.

\begin{figure}
	\centering
	\includegraphics[width=0.4\textwidth]{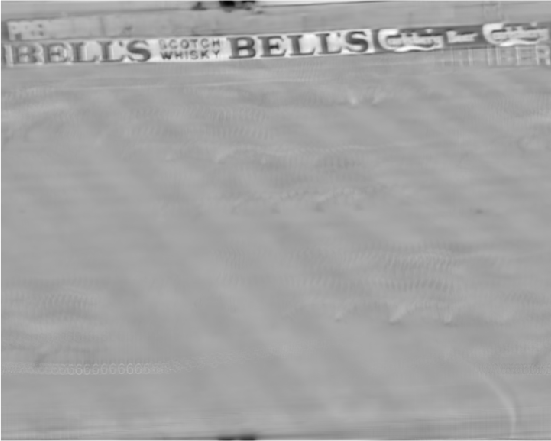}\\[0.1cm]
	\includegraphics[width=0.4\textwidth]{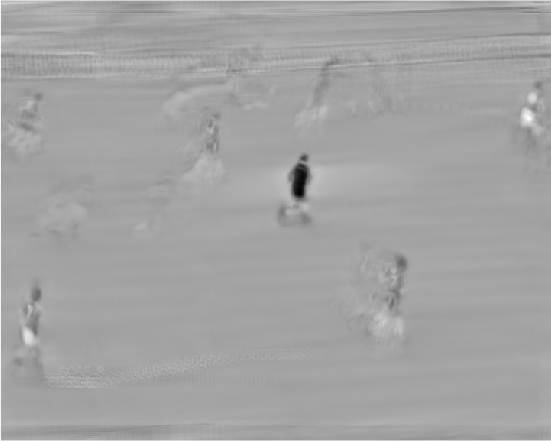}\\[0.1cm]
	\includegraphics[width=0.4\textwidth]{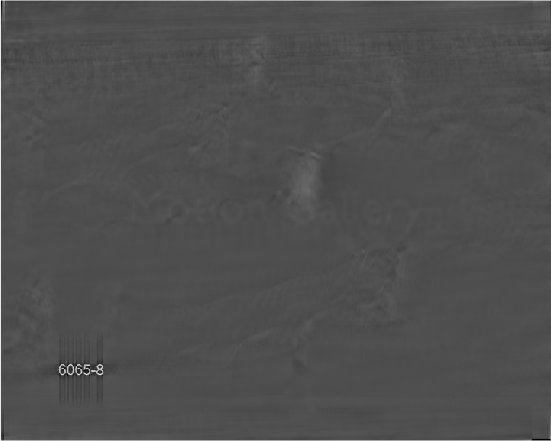}
	\caption{Tracked objects in soccer clip: advertising banners, referee and time stamp.}
	\label{fig:video2objects}
\end{figure}

However, we were still able to extract information using our method with $L=10$. The first extracted object is mainly background. Objects two, three and four form the advertising banners. The next object is the referee together with some artifacts from players who move at same speed as the referee. Finally, as 10th object we were able to extract the time stamp in the lower left corner, i.e., the only object that does not move during the video. Figure \ref{fig:video2objects} shows the extracted objects.

\section{Conclusion and future work}

We presented a generalization of the low-rank approximation, which allows to individually shift the column of rank-$1$ matrices. This model was designed to represent objects that move through the data. This holds in applications such as seismic or ultrasonic image processing as well as video processing. Basic properties of the occurring shift operators and shifted rank-$1$ matrices were stated. With the help of these properties we were able to develop an efficient algorithm to find a shifted rank-$1$ approximation of given input. The new approach comes with some disadvantages that we want to overcome in future works. The sparse approximation is slightly worse than comparable methods. Also, objects moving in and out of a video, changing in shape or covering other objects may cause problems. However, big advantages of this model are the extracted parameters. Their importance for applications has been demonstrated in several examples including approximation, segmentation and object tracking. We compared the results to other methods such as wavelet transform or shift-invariant dictionary learning. It has been shown that the obtained parameters can directly be used to extract crucial information. Since we make extensive use of properties of the Fourier transform, the approximation algorithm has a complexity of only $O(LMN\log M)$. This also holds for the matrix vector multiplication of a shifted rank-$1$ matrix that can be performed in $O(M\log M+N)$.

We plan to extend the method in future works to overcome some of the drawbacks of the model. Integrating conditions such as compact support of one or both of the singular vectors or smoothness of the shift vector are of special interest. Furthermore, besides shifting the columns, other operators, e.g., scaling, may be applied.

\section{Acknowledgments}
The work was supported in part by National Key Research and Development Program of China under Grant 2017YFB0202900, and NSFC under Grant 41625017 and Grant 41374121.

\bibliographystyle{IEEEtran}
\bibliography{references}

\end{document}